\documentclass[a4paper,11pt,twoside]{article}

\usepackage{amssymb}
\usepackage[T1]{fontenc}
\usepackage[utf8]{inputenc}
\usepackage{mathrsfs}
\usepackage{amsfonts,amsmath,lmodern,fancyhdr,lastpage,graphicx,nccfoots,caption}
\usepackage{afterpage}

\usepackage{rlepsf}

\usepackage{graphicx}
\usepackage[all]{xy}
\usepackage{amscd}
\usepackage{comment}

\newtheorem{Theorem}{Theorem}[section]

\newtheorem{Definition}[Theorem]{Definition}

\newtheorem{Proposition}[Theorem]{Proposition}
\newtheorem{Corollary}[Theorem]{Corollary}

\newtheorem{Remark}[Theorem]{Remark}

\newenvironment{Proof}[1][Proof]{\textbf{#1.} }{\ \rule{0.5em}{0.5em}}

\def \R {\mathbb{R}}

\def \be {\begin{equation}}

\def \ee {\end{equation}}

\usepackage{hyperref}

\vfuzz \hfuzz
\newcounter{a}
\ifodd\thea\else\stepcounter{a}\fi

\fancypagestyle{plain}{%
\fancyhf{}
}
\begin{document}
%\label{pp}
\thispagestyle{plain}

%%%%%%%%%%%%%%%%%%%%
% Ttulo do artigo
\begin{center}
\Large
\textsc{Invariants and TQFT's for cut cellular surfaces from finite groups}
\end{center}

%%%%%%%%%%%%%%%%%%%%
% Primeiros autores, com endereo comum

\begin{center}
\textit{Diogo Bragan\c ca}  \smallskip \\
\begin{tabular}{l}
\small Department of Physics, \\ 
\small  Instituto Superior T\'ecnico, Universidade de Lisboa \\
\small Av. Rovisco Pais, 1                    \\
\small 1049-001 Lisboa, Portugal       \\
\small e-mail: \texttt{diogo.braganca@hotmail.com}   \\ 
\end{tabular}
\end{center}

%%%%%%%%%%%%%%%%%%%%
% Outro autor, com endereo diferente

\begin{center}
\textit{Roger Picken} \smallskip \\
\begin{tabular}{l}
\small Center for Mathematical Analysis, Geometry and Dynamical Systems, \\
\small Department of Mathematics, \\
\small  Instituto Superior T\'ecnico, Universidade de Lisboa                \\
\small Av. Rovisco Pais, 1                    \\
\small 1049-001 Lisboa, Portugal      \\
\small e-mail: \texttt{roger.picken@tecnico.ulisboa.pt}
\end{tabular}
\end{center}

%%%%%%%%%%%%%%%%%%%%
% Resumo (em portugus)

\medskip

%%%%%%%%%%%%%%%%%%%%
% Abstract (em ingl\^es)  

\noindent
\textbf{Abstract} We introduce the notion of a cut cellular surface (CCS), being a surface with boundary, which is cut in a specified way to be represented in the plane, and is composed of 0-, 1- and 2-cells. We obtain invariants of CCS's under Pachner-like moves on the cellular structure, by counting colourings of the 1-cells with elements of a finite group, subject to a ``flatness'' condition for each 2-cell. These invariants are also described in a TQFT setting, which is not the same as the usual 2-dimensional TQFT framework. We study the properties of functions which arise in this context, associated to the disk, the cylinder and the pants surface, and derive general properties of these functions from topology, including properties which come from invariance under the Hatcher-Thurston moves on pants decompositions.

\medskip

%%%%%%%%%%%%%%%%%%%%
% keywords (vers\~ao inglesa da lista de palavras-chave)

\noindent\textbf{keywords:} Cut cellular surface, TQFT, finite group, pants decomposition

%%%%%%%%%%%%%%%%%%%%
%%%%%%%%%%%%%%%%%%%%
% Primeira seco

\newpage

\section{Introduction}

The purpose of this article is to describe a class of invariants and TQFT's for a particular type of surfaces, that we call cut cellular surfaces. Since the notion of TQFT (Topological Quantum Field Theory), may not be familiar to all readers, we start, in Section 2, with a brief introduction to this subject, based around a simple example. 

In Section 3 we introduce our central concept, the definition of a cut cellular surface (CCS). This is a surface with boundary, which is cut in a specified way to be represented in the plane (like the well-known rectangle with opposite edges identified representing the torus), and which is composed of 0-, 1- and 2-cells, generalizing the familiar notion of a triangulated surface. We provide examples and define two moves on the cell structure which give equivalent planar representations. In the special case of triangulated manifolds without boundary we show that these moves generate the Pachner moves on triangulations. 

In Section 4, we define invariants of CCS's which come from assignments of elements of a finite group $G$ to the 1-cells of the CCS. The assignments, called $G$-colourings, are subject to a condition (``flatness'') and the invariant counts the number of valid $G$-colourings with a suitable normalization factor. These expressions are shown to be invariant under the two moves on the cellular structure. We calculate the invariant for some elementary examples. 

In Section 5 we describe how the invariant behaves when gluing two CCS's together along a common boundary component, and use this to get a TQFT for these surfaces. We note that these TQFT's are not the same as the 2D TQFT's which were classified by Abrams in terms of Frobenius algebras \cite{ab}.

In Section 6 we analyse general properties of the invariants, such as their behaviour under certain symmetry operations on CCS's. We study properties of the invariants assigned to three elementary surfaces, the disk $D$, the cylinder $C$ and the pants surface $P$, and give a sample calculation of how a topological equivalence implies a relation for  these invariants. We conclude by giving two properties of the invariant assigned to the pants surface, which correspond to the Hatcher-Thurston moves \cite{ht} on pants decompositions of a surface.   

As we point out in our conclusions, in Section 7, this article paves the way for an analogous construction, using colourings of both 1- and 2-cells of a CCS, with elements of a finite 2-group.

\vskip 0.3cm

\section{A brief introduction to TQFT}

The notion of TQFT (Topological Quantum Field Theory), which motivates our approach, may not be familiar to all readers, so we start with a brief introduction to this subject. In essence, TQFT is a way of obtaining invariants of a topological or geometrical nature, based on the formalism used by physicists in describing quantum theory.

Suppose we have a class of manifolds with boundary, and we specify for each $M$ in this class, an ``in'' part and ``out'' part of its boundary. Such a manifold represents the evolution of a physical system from an initial to a final configuration. The idea is that we will only take into account basic topological or geometrical features of $M$, i.e. we do not distinguish between $M$'s with equivalent features. 

\begin{figure}[htbp] 
\centerline{\relabelbox 
\epsfysize 1.5cm
\epsfbox{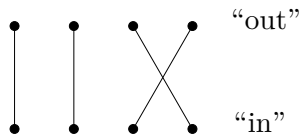}
\relabel{i}{{\rm``in''}}
\relabel{o}{{\rm``out''}}
\endrelabelbox}
\caption{\label{fig:permutn} A manifold $M$ with 4 strands}
\end{figure}

As an example, let $M$ be a union of 1-dimensional manifolds (``strands'') connecting $n$ ordered points (the ``in'' boundary) to $n$ ordered points (the ``out'' boundary), but not necessarily in the same order. See Figure \ref{fig:permutn} below representing such a manifold with 4 strands. 

Note that the two strands connecting the 3rd and 4th points do not actually intersect, despite the appearance in the planar representation. 

\begin{figure}[htbp] 
\centerline{\relabelbox 
\epsfysize 2.3cm
\epsfbox{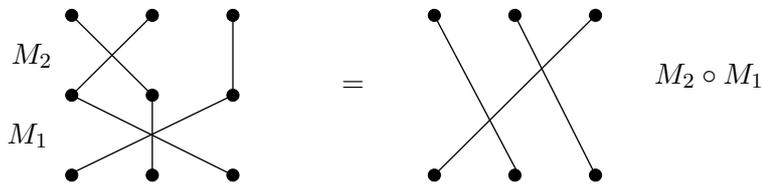}
\relabel{i}{$M_1$}
\relabel{o}{$M_2$}
%\relabel{k}{$M_2\circ M_1$}
\adjustrelabel <-2pt,3pt> {k}{$M_2\circ M_1$}
\relabel{l}{$=$}
\endrelabelbox}
\caption{\label{fig:permutn-comp} Composition of manifolds }
\end{figure}

For two such manifolds, $M_1$ and $M_2$, with the same number of strands $n$, we can perform a composition or concatenation, denoted $M_2\circ M_1$, as in Figure \ref{fig:permutn-comp}, where $n=3$. The two representations of this composition in the Figure are equivalent, since we are only concerned with the topology, which in this case means no more than what the endpoints of each strand are.

Thus for fixed $n$, we clearly have a correspondence between equivalence classes of manifolds $M$, and elements of the permutation group $P_n$. Composition of manifolds corresponds to multiplication in this group.

A TQFT is an algebraic representation of the class of manifolds $M$ in the following sense. To each type of ``in'' or ``out'' boundary, denoted $\Sigma$, we assign a finite-dimensional vector space 
$V_\Sigma$ (over some fixed field $k$). In the language of quantum theory, this vector space represents the Hilbert space of states for a physical system. To the manifold $M$ itself, we assign a linear transformation $Z_M: V_{\Sigma_i}\rightarrow V_{\Sigma_o}$. In the context of quantum theory, $Z_M$ corresponds to the evolution operator of the physical system.

In our example, if we denote by $V_n$ the vector space assigned to $n$ points, a manifold $M$ with $n$ strands gets assigned to it a linear transformation $Z_M:V_n \rightarrow V_n$.

The assignments $\Sigma\mapsto V_\Sigma$ and $M\mapsto Z_M$ need to satisfy a number of natural properties in order to constitute a TQFT. We will just give the main ones, which are sufficient for the understanding of the present article, and refer the reader to the article by Atiyah \cite{at}, where the axioms of TQFT were first introduced. 

First, when two manifolds $M_1$ and $M_2$ are composed, we have the property:
\be
Z_{M_2\circ M_1} = Z_{M_2}\circ Z_{M_1}
\label{eq:comp}
\ee
i.e. the composition of manifolds is represented by the composition of the corresponding linear transformations. In the case of our example, this means that the linear maps $Z_M$, for manifolds $M$ with $n$ strands, constitute a representation of $P_n$ on the vector space $V_n$.

An immediate consequence of property (\ref{eq:comp}) follows from applying it to $M=\Sigma \times I$, where $I$ is a closed interval in $\R$, and $M$ has ``in''  and ``out'' boundary both equal to $\Sigma$. Since $(\Sigma \times I) \circ (\Sigma \times I)$ is equivalent to $\Sigma \times I$, we get the relation:
\be
Z_{\Sigma \times I}=  Z_{(\Sigma \times I) \circ (\Sigma \times I)} = Z_{\Sigma \times I} \circ  Z_{\Sigma \times I},
\label{eq:idemp}
\ee
i.e. $Z_{\Sigma \times I}$ is an idempotent. For this reason, in  many TQFT's, $Z_{\Sigma \times I}$ is taken to be ${\rm id}_{V_\Sigma}$, the identity map on $V_\Sigma$, and this is the choice we will make in our example. Thus we will assign to the manifold consisting of $n$ vertical strands (like the two strands on the left in Figure \ref{fig:permutn}) the identity transformation on $V_n$.

The second main property of a TQFT refers to the situation where either $\Sigma$ or $M$ are disjoint unions of manifolds. The disjoint union is the union of disjoint sets, e.g. $\Sigma$ consisting of 5 points is the disjoint union of $\Sigma_1$ consisting of 2 points and $\Sigma_2$ consisting of 3 points. Likewise the manifold $M$ of Figure \ref{fig:permutn} is the disjoint union of $M_1$, consisting of a pair of parallel strands, and $M_2$, consisting of a pair of crossing strands. We will write this as: $\Sigma=\Sigma_1\sqcup \Sigma_2$ or $M=M_1\sqcup M_2$.

For $\Sigma=\Sigma_1\sqcup \Sigma_2$ we have the property
\be
V_{\Sigma_1\sqcup \Sigma_2}=V_{\Sigma_1} \otimes V_{\Sigma_2},
\label{eq:VSigma}
\ee
where the right hand side denotes the tensor product of vector spaces. Since the vector spaces $V_{\Sigma_1}$ and $V_{\Sigma_2}$ are finite-dimensional, given bases $\left\{e_i\right\}_{i=1,\dots,p}$ and 
$\left\{f_j\right\}_{j=1,\dots,q}$ of $V_{\Sigma_1}$ and $V_{\Sigma_2}$ respectively, we have a basis of $V_{\Sigma_1} \otimes V_{\Sigma_2}$ consisting of $pq$ elements of the form $e_i\otimes f_j$.

Likewise, when $M$ is a disjoint union $M=M_1\sqcup M_2$, we have the property
\be
Z_{M_1\sqcup M_2} = Z_{M_1}\otimes Z_{M_2},
\label{eq:ZM}
\ee
where, on the right hand side, the tensor product of linear transformations is defined by:
$$
(Z_{M_1}\otimes Z_{M_2}) (v\otimes w) = Z_{M_1}(v)\otimes Z_{M_2} (w).
$$
Note that, since $M$ is a disjoint union, the ``in''  and ``out'' boundaries of $M$ are disjoint unions too.

In our example, the property (\ref{eq:VSigma}) means that we only need to specify $V_\Sigma=V$, when $\Sigma$ is a single point, since for $\Sigma$ equal to $n$ points we have $V_n=V\otimes V \cdots \otimes V$ ($n$ times), or  more succinctly, $V_n = V^{\otimes n}$. Likewise, for $M=M_1\sqcup M_2$, the linear map $Z_M$ is constrained by the property (\ref{eq:ZM}). 

We can summarise all the above by saying that a TQFT for our example consists of a compatible collection of representations of $P_n$ on $V_n = V^{\otimes n}$, for all $n$. An example of such a collection is given as follows: let $M$ be an $n$-strand manifold which connects point $i$ of the ``in '' boundary to point $\sigma(i)$ of the ``out '' boundary, for $i=1,\dots , n$, where $\sigma\in P_n$. Then we set:
$$
Z_M(v_1\otimes \cdots \otimes v_n) = v_{\sigma(1)}\otimes v_{\sigma(1)} \cdots \otimes v_{\sigma(n)}
$$ 
for any elements $v_1, \dots, v_n\in V$, extended by $k$-linearity to general elements of $V^{\otimes n}$. This clearly satisfies the conditions (\ref{eq:comp}), (\ref{eq:idemp}) and (\ref{eq:ZM}) for $Z_M$, and (\ref{eq:VSigma}) has already been incorporated by having 
$V_n=V^{\otimes n}$.

We haven't touched here on some other general properties of a TQFT, such as those which relate to the cases when $\Sigma$ or $M$ is the empty set,
or the behaviour under change of orientation when $\Sigma$ and $M$ are oriented manifolds. These properties are irrelevant for the discussion of our example, but come into play when considering more general types of manifold, such as an oriented $\cup$-shaped manifold, with ``in'' boundary the empty set and ``out'' boundary two oppositely-oriented points. The empty manifolds $\Sigma$ and $M$ are also needed in order to complete the formal mathematical structure surrounding equations 
(\ref{eq:VSigma}) and (\ref{eq:ZM}). We refer again to \cite{at} for an in-depth presentation of the general properties of a TQFT.

Finally, we will say a brief word about notation. When we wish to describe a linear transformation $Z_M: V_{\Sigma_i}\rightarrow V_{\Sigma_o}$ in concrete terms, we may introduce a basis 
$\left\{e_i\right\}_{i=1,\dots,n}$ of $V_{\Sigma_i}$ and a basis $\left\{f_j\right\}_{j=1,\dots,m}$ of $V_{\Sigma_o}$. Then $Z_M$ is represented by an $m\times n$ matrix 
$[c_{ji}]$, where $Z_M(e_i) =\sum_{j=1}^m c_{ji}f_j$. We will be using the suggestive physicists' notation for the matrix elements $c_{ji}$, namely:
$$
c_{ji}=  \left \langle f_j\left | Z_M \right | e_i  \right \rangle.
$$

\section{Cut cellular surfaces}
\label{sec:a}

\begin{figure}[htbp]
\centering
\includegraphics[width=10cm]{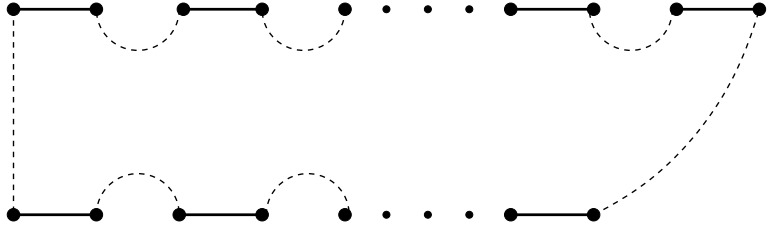}
\caption{General appearance of a cut cellular surface (CCS)}
\label{fig:ccs}
\end{figure}

We will be considering surfaces with boundary, which are cut in a specified way to be represented in the plane (like the well-known rectangle with opposite edges identified representing the torus), and which are composed of 0-, 1- and 2-cells, generalizing the familiar notion of a triangulated surface. We recommend using Figure \ref{fig:ccs} to accompany the following definition.

%%%%%

\begin{Definition}
A cut cellular surface (CCS) is an orientable 2-manifold $M$ with boundary, endowed with a finite cell-structure, such that
\begin{itemize}
\item[a)] Each boundary component of $M$ consists of a single 0-cell and a single 1-cell.
\item[b)] $M$ has a specified planar representation, obtained by cutting $M$ along 1-cells in such a way as to obtain a simply connected region in the plane, bounded by ``external'' 1-cells, namely the boundary 1-cells and the 1-cells along which $M$ was cut (the latter appear twice in the boundary of the planar region). The cut 1-cells are labelled and given an orientation to make explicit how they are identified in $M$.
\item[c)] The planar representation has the schematic structure shown in Fig. \ref{fig:ccs}: the boundary components, represented by solid lines, lie either along the bottom or the top edge of the planar representation. Those along the bottom edge are called ``in'' boundary components, those along the top edge are called ``out'' boundary components. When there are no ``in/out'' boundary components, the bottom/top edge contains a single 0-cell. The dotted lines on the left and right, and the dotted lines between boundary components along the bottom and top edge, each represent one or more cut 1-cells, separated by 0-cells when there are more than one of them.
\item[d)] The simply connected planar region is made up of one or more 2-cells, separated by ``internal'' 1-cells and 0-cells.

\end{itemize}
\label{def:ccs}
\end{Definition}

\begin{figure}[htbp] 
\centerline{\relabelbox 
\epsfysize 2.3cm
\epsfbox{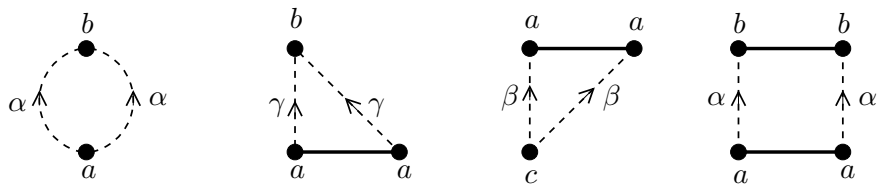}
\relabel{a}{$a$}
\relabel{b}{$b$}
\relabel{c}{$a$}
\adjustrelabel <-1pt,0pt> {d}{$a$}
\relabel{e}{$b$}
\relabel{f}{$c$}
\relabel{g}{$a$}
\adjustrelabel <0pt,1pt> {h}{$a$}
\relabel{i}{$a$}
\adjustrelabel <-1pt,0pt> {j}{$a$}
\relabel{k}{$b$}
\relabel{l}{$b$}
\relabel{X}{$\alpha$}
\relabel{Y}{$\alpha$}
\relabel{Z}{$\gamma$}
\relabel{W}{$\gamma$}
\adjustrelabel <-1pt,0pt> {T}{$\beta$}
\relabel{U}{$\beta$}
\relabel{R}{$\alpha$}
\relabel{S}{$\alpha$}
\endrelabelbox}
\caption{\label{fig:SDC} The sphere $S$, disc $D$ and cylinder $C$ represented as CCS's}
\end{figure}

\begin{figure}[htbp] 
\centerline{\relabelbox 
\epsfysize 3cm
\epsfbox{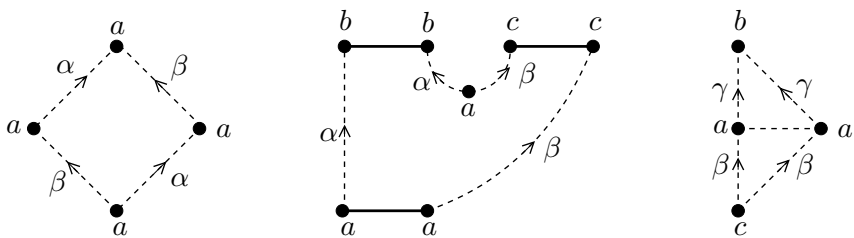}
\relabel{a}{$a$}
\adjustrelabel <-1pt,0pt> {b}{$a$}
\adjustrelabel <-1pt,0pt> {p}{$a$}
\relabel{q}{$a$}
\relabel{D}{$\alpha$}
\adjustrelabel <-1pt,1pt> {B}{$\alpha$}
\adjustrelabel <-2pt,-2pt> {A}{$\beta$}
\relabel{C}{$\beta$}
\relabel{i}{$a$}
\adjustrelabel <-2pt,0pt> {j}{$a$}
\adjustrelabel <-2pt,0pt> {k}{$b$}
\adjustrelabel <-2pt,1pt> {l}{$b$}
\relabel{m}{$a$}
\relabel{c}{$c$}
\relabel{d}{$c$}
\relabel{R}{$\alpha$}
\relabel{S}{$\beta$}
\adjustrelabel <-1pt,-1pt> {X}{$\alpha$}
\relabel{Y}{$\beta$}
\relabel{e}{$b$}
\relabel{f}{$c$}
\relabel{g}{$a$}
\relabel{h}{$a$}
\adjustrelabel <-2pt,0pt> {T}{$\beta$}
\relabel{U}{$\beta$}
\adjustrelabel <-2pt,0pt> {Z}{$\gamma$}
\relabel{W}{$\gamma$}
\endrelabelbox}
\caption{\label{fig:TPS} The torus $T$, pants surface $P$, and sphere $S$, now as two glued discs}
\end{figure}

\begin{figure}[htbp] 
\centerline{\relabelbox 
\epsfysize 4cm
\epsfbox{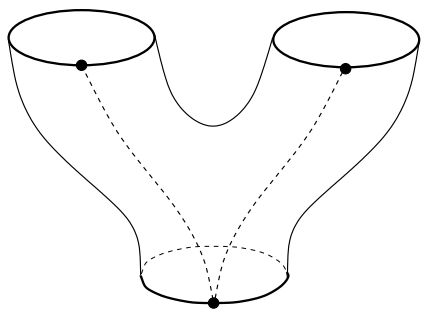}
\endrelabelbox}
\caption{\label{fig:3D-pants} The pants surface and its cuts}
\end{figure}

To fix ideas we give some simple examples. 
\begin{itemize}
\item The example with the least number of 1-cells is the sphere $S$, represented on the left in Fig. \ref{fig:SDC} with two 0-cells, one cut 1-cell and one 2-cell. 

\item The disc $D$ can be represented as a triangle in two ways, depending on whether the boundary is taken to be ``in'' or ``out'' (the two middle CCS's in Fig. \ref{fig:SDC}). 

\item A representation of the cylinder $C$ is shown on the right in Fig. \ref{fig:SDC}.

\item The well-known representation of the torus $T$ as a square with opposite edges identified is shown as a CCS on the left in Figure \ref{fig:TPS}.

\item We show a representation of the pants surface $P$, choosing (for sartorial reasons!) one ``in'' boundary component and two ``out'' boundary components (in the middle of Fig. \ref{fig:TPS}). We illustrate the cuts made in Figure \ref{fig:3D-pants}.

\item If we glue together two discs to make a sphere, as shown on the right in Fig. \ref{fig:TPS}, we get an example with more than one 2-cell. This leads us to study moves between different planar representations of the same surface, such as the two different representations of the sphere $S$ in Figs \ref{fig:SDC} and \ref{fig:TPS}.
\end{itemize}

{\bf Moves on CCS's.} By analogy with the Pachner moves on triangulated manifolds, we introduce moves for passing between different planar representations of the same surface. There are two types of move. 
\vskip 0.3cm
 
\begin{figure}[htbp] 
\centerline{\relabelbox 
\epsfxsize 8cm
\epsfbox{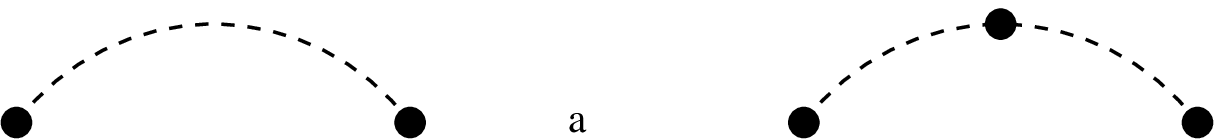}
\relabel{a}{$\longleftrightarrow$}
\endrelabelbox}
\caption{\label{fig:move1} Move I}
\end{figure}

{\bf Move I:} Introducing a 0-cell into a 1-cell, thereby dividing it into two 1-cells, or conversely removing a 0-cell separating two 1-cells, to combine them into a single 1-cell (Figure \ref{fig:move1}).

\begin{figure}[htbp] 
\centerline{\relabelbox 
\epsfxsize 13cm
\epsfbox{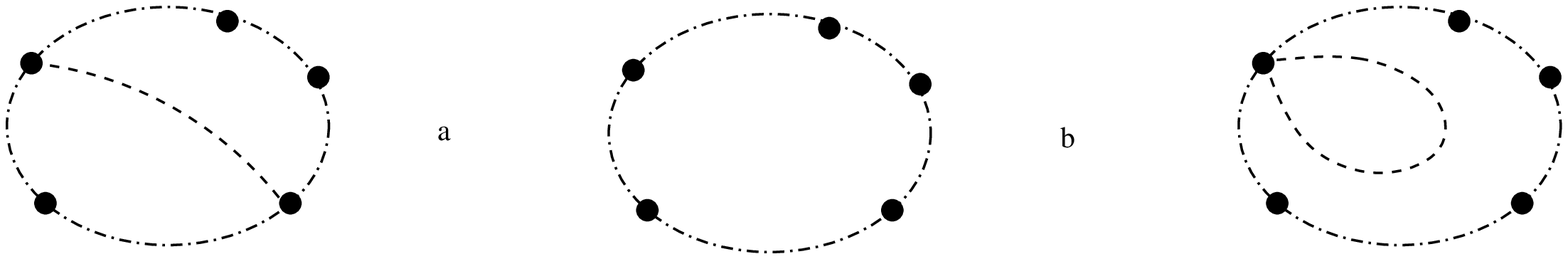}
\relabel{a}{$\longleftrightarrow$}
\relabel{b}{$\longleftrightarrow$}
\endrelabelbox}
\caption{\label{fig:move2} Move II}
\end{figure}

{\bf Move II:} Introducing a 1-cell into a 2-cell, thereby dividing it into two 2-cells, or conversely removing a 1-cell separating two 2-cells, to combine them into a single 2-cell (Figure \ref{fig:move2}).

\begin{Remark}
Regarding move I, it applies only to the 1-cells which are not boundary components of $M$, since we have made it a rule that the boundary components consist of a single 1-cell and 0-cell, see Def. \ref{def:ccs} a). When this move is applied to a cut 1-cell, we must ensure that the 0-cell is introduced into both copies of the cut 1-cell, or removed from between both copies of the cut 1-cell, in the planar representation. 

Regarding Move II, we have used dots-and-dashes lines for the 1-cells around the perimeter of the 2-cell to indicate that they can be either boundary or non-boundary 1-cells. Note that the 1-cell which is introduced or removed, connects two 0-cells in the boundary of the 2-cell in the planar representation, and that these may be the same 0-cell, as shown on the right in Figure \ref{fig:move2}.  
\end{Remark}

Clearly a triangulated surface without boundary, endowed with a planar representation, is an example of a CCS. As is well-known, different triangulations of the same surface are related by the Pachner moves \cite{pachner}. These moves are generated by the moves on CCS's, as we show in Figure \ref{fig:pachner} below.

\begin{figure}[htbp]
\centering
\includegraphics[width=11cm]{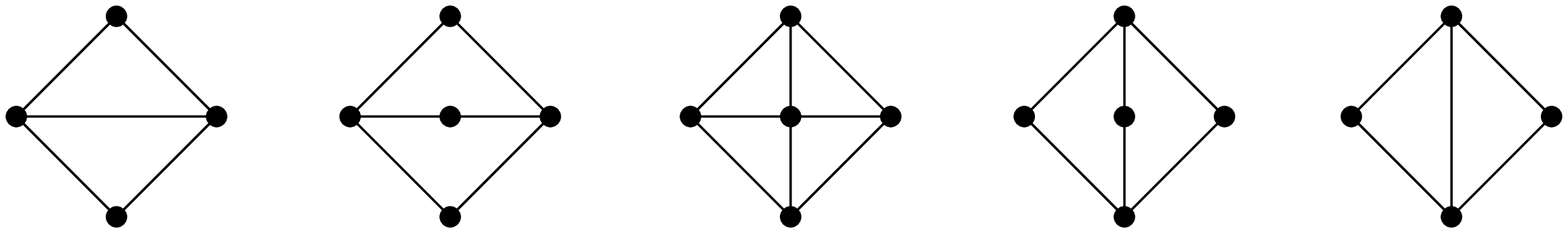}
\vskip 0.9cm
\includegraphics[width=12cm]{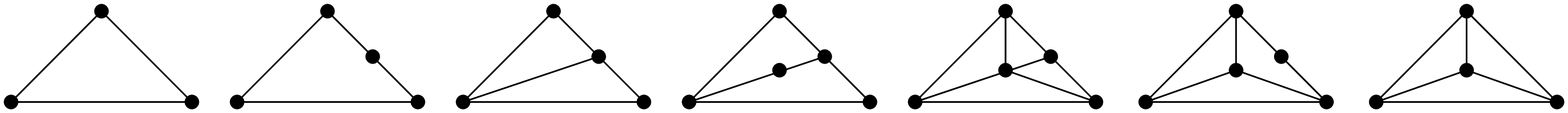}
\caption{The Pachner moves as sequences of moves on CCS's}
\label{fig:pachner}
\end{figure}

\section{Invariants for CCS's from finite group colourings}
\label{sec:b}

Fix a finite group $G$. Given a CCS, $M$, we fix orientations on the 1-cells of $M$, specified as follows with respect to the planar representation:
\begin{itemize}
\item the boundary 1-cells are oriented from left to right
\item the cut 1-cells are oriented as chosen in Definition \ref{def:ccs} b)
\item the internal 1-cells are oriented arbitrarily.
\end{itemize}

\begin{Definition} A $G$-colouring of $M$ is an assignment of an element $g_i\in G$ to each 1-cell labelled $i$, such that, for each 2-cell in the planar representation, the following (flatness) condition holds:
\begin{itemize}
\item if the 1-cells of the boundary of the 2-cell are labelled $i_1, \dots i_k$, ordered in the anticlockwise direction, then
\begin{equation}
\prod_{j=1}^k g_{i_j}^{(-1)} =1
\end{equation} 
where the factor is $g_{i_j}$ or $g_{i_j}^{-1}$, depending on whether or not the 1-cell $i_j$ is oriented compatibly with the positive orientation of the 2-cell.
\end{itemize}
\end{Definition}

\begin{figure}[htbp] 
\centerline{\relabelbox 
\epsfxsize 3cm
\epsfbox{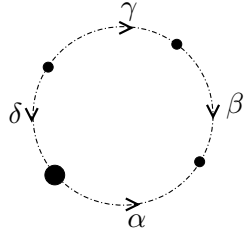}
\relabel{a}{$\alpha$}
\relabel{b}{$\beta$}
\relabel{c}{$\gamma$}
\relabel{d}{$\delta$}
\endrelabelbox}
\caption{The flatness condition here is $g_\alpha g_\beta^{-1}g_\gamma^{-1}g_\delta =1$.}
\label{fig:flatness}
\end{figure}

See Figure \ref{fig:flatness} for an example of the flatness condition. We have again used dots-and-dashes lines for the 1-cells to indicate that they can be either boundary or non-boundary 1-cells.

\begin{Remark}
Clearly the flatness condition does not depend on the choice of 0-cell from which we start ordering the 1-cells, e.g. in the example of Figure \ref{fig:flatness} where the starting 0-cell has been shown thickened, the condition $g_\alpha g_\beta^{-1}g_\gamma^{-1}g_\delta =1$ is equivalent to the condition 
$g_\beta^{-1}g_\gamma^{-1}g_\delta g_\alpha  =1$, which corresponds to starting instead at the 0-cell which is the endpoint of $\alpha$.
The term flatness condition is motivated by the links with gauge theory, where a flat $G$-connection gives rise to $G$-valued transports along oriented 1-cells satisfying the flatness condition around each 2-cell.
\end{Remark}

We can define an invariant of CCS's, using $G$-colourings. Choose elements of $G$, $g_1,\dots , g_n$, for the colouring of the ``in'' boundary components, and 
$h_1,\dots , h_m$, for the colouring of the ``out'' boundary components, ordering the boundary components from left to right in the planar representation. Let  $|G|$ denote the number of elements of the finite group $G$, and $v$ denote the number
of internal vertices, i.e. 0-cells, of $M$. Let ${\cal C}={\cal C}(g_1,\dots , g_n; h_1,\dots , h_m)$ denote the set of all $G$-colourings of $M$ which have the given assignments on the boundary components.
Then we define:
\begin{equation}
\left \langle h_1,\dots , h_m \left | Z_M \right |g_1,\dots , g_n   \right \rangle
:= \frac{1}{|G|^{\frac{m+n}{2}+v}} \, \# {\cal C}(g_1,\dots , g_n; h_1,\dots , h_m).
\label{def:Zgrp}
\end{equation}
If there are no ``in'' or ``out'' components we write the invariant as  $\left \langle \dots  \left | Z_M \right |\emptyset  \right \rangle$ or $\left \langle \emptyset \left | Z_M \right | \dots    \right \rangle$. 

\begin{Remark}
Equivalently we can write the invariant as a state-sum:
$$
\left \langle h_1,\dots , h_m \left | Z_M \right |g_1,\dots , g_n   \right \rangle := \frac{1}{|G|^{\frac{m+n}{2}+v}} \,\sum_{c\in {\cal C}} 1 .
$$
\end{Remark}

We now discuss in what sense these are invariants. First of all, we have:
\begin{Proposition}
The invariants $\left \langle h_1,\dots , h_m \left | Z_M \right |g_1,\dots , g_n   \right \rangle$ are unchanged under changes of orientation of the cut 1-cells or internal 1-cells.
\label{prop:orientn}
\end{Proposition}
\begin{Proof}
Clearly the number of internal vertices is unchanged, and there is a bijection between the respective sets of colourings, given by replacing the element $g$ assigned to any cut 1-cell or internal 1-cell by $g^{-1}$, when its orientation is reversed.
\end{Proof}
\vskip 0.1cm

\begin{figure}[htbp] 
\centerline{\relabelbox 
\epsfxsize 8cm
\epsfbox{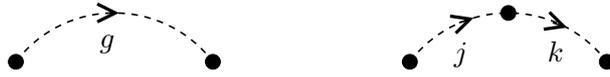}
\relabel{g}{$g$}
\relabel{j}{$j$}
\relabel{k}{$k$}
\endrelabelbox}
\caption{Colourings of $M$ and $M'$ for Move I}
\label{fig:pfmv1}
\end{figure}

\begin{figure}[htbp] 
\centerline{\relabelbox 
\epsfxsize 11.5cm
\epsfbox{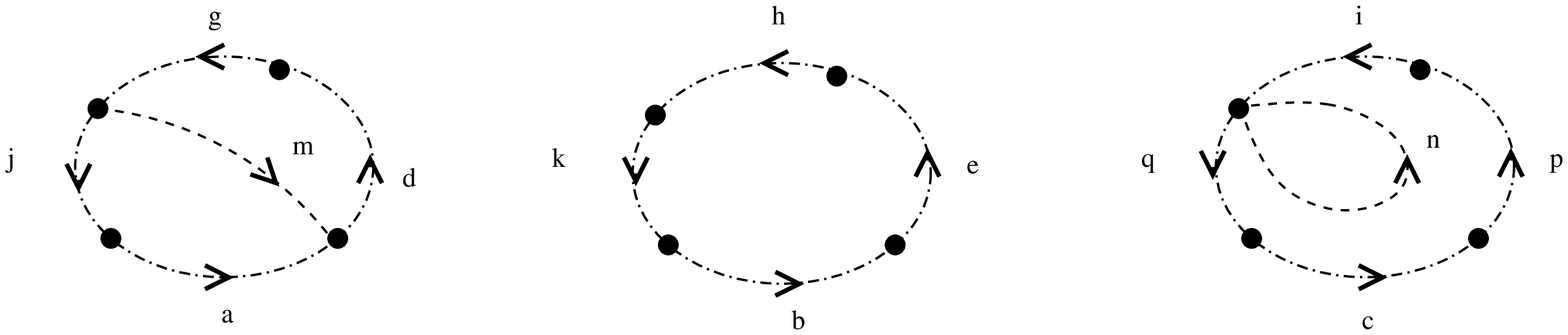}
\relabel{a}{$k_1$}
\relabel{b}{$k_1$}
\relabel{c}{$k_1$}
\relabel{d}{$k_2$}
\relabel{e}{$k_2$}
\relabel{p}{$k_2$}
\relabel{g}{$k_3$}
\relabel{h}{$k_3$}
\relabel{i}{$k_3$}
\relabel{j}{$k_4$}
\relabel{k}{$k_4$}
\relabel{q}{$k_4$}
\adjustrelabel <-1pt,0pt> {m}{$g$}
\adjustrelabel <-1pt,-1pt> {n}{$h$}
\endrelabelbox}
\caption{Colourings of $M$ and $M'$ for Move II}
\label{fig:pfmv2}
\end{figure}

More importantly we have:
\begin{Theorem}
The invariants $\left \langle h_1,\dots , h_m \left | Z_M \right |g_1,\dots , g_n   \right \rangle$ are unchanged under moves I and II.
\end{Theorem}
\begin{Proof} Suppose $M$ and $M'$ are related by a move I. Fix a $G$-colouring for $M$ that assigns $g$ to the 1-cell displayed on the left in Figure \ref{fig:pfmv1}. Keeping the assignments of all other 1-cells the same, for $M'$ on the right there are $|G|$ compatible $G$-colourings, since we can choose one assignment, e.g. $j$, freely in $G$ and the other assignment $k$ is then determined (for the orientations as shown in Figure \ref{fig:pfmv1}, we have $g=jk$, i.e. $k=j^{-1}g$). Since $M'$ has an extra internal vertex compared to $M$ the invariants (\ref{def:Zgrp}) are the same for $M$ and $M'$.

We prove the invariance under move II for an obviously representative case, shown in Fig. \ref{fig:pfmv2}. For each $G$-colouring of the middle 2-cell, there is precisely one compatible $G$-colouring for the 2-cells on the left or the right (fixing the assignments on all 1-cells not displayed). For the orientations in the figure we have $g=k_1^{-1}k_4^{-1}=k_2k_3$ and $h=1$. Since the number of internal vertices is the same in all three cases, so are the invariants (\ref{def:Zgrp}).  
\end{Proof}

\vskip 0.3cm

\begin{figure}[htbp] 
\centerline{\relabelbox 
\epsfysize 2.8cm
\epsfbox{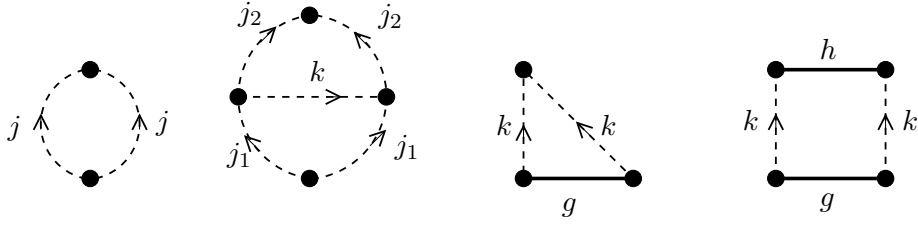}
\relabel{X}{$j$}
\relabel{Y}{$j$}
\relabel{A}{$j_1$}
\adjustrelabel <-1pt,1pt> {B}{$j_1$}
\relabel{C}{$j_2$}
\adjustrelabel <-1pt,-1pt> {D}{$j_2$}
\relabel{E}{$k$}
\relabel{Z}{$k$}
\relabel{W}{$k$}
\relabel{R}{$k$}
\relabel{S}{$k$}
\relabel{T}{$g$}
\relabel{U}{$g$}
\relabel{h}{$h$}
\endrelabelbox}
\caption{\label{fig:inv-exp1} $G$-colourings for the sphere, disk and cylinder}
\end{figure}

\begin{figure}[htbp] 
\centerline{\relabelbox 
\epsfysize 3.5cm
\epsfbox{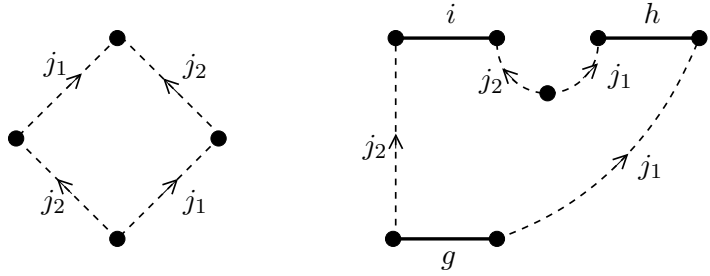}
\relabel{g}{$g$}
\relabel{i}{$i$}
\relabel{h}{$h$}
\relabel{D}{$j_1$}
\adjustrelabel <-1pt,1pt> {B}{$j_1$}
\relabel{A}{$j_2$}
\relabel{C}{$j_2$}
\relabel{Y}{$j_1$}
\relabel{S}{$j_1$}
\adjustrelabel <-1pt,0pt> {X}{$j_2$}
\adjustrelabel <-1pt,0pt> {R}{$j_2$}
\endrelabelbox}
\caption{\label{fig:inv-exp2} $G$-colourings for the torus and the pants surface}
\end{figure}

Next we calculate the invariant for some simple examples, referring to Figures  \ref{fig:inv-exp1} and \ref{fig:inv-exp2}.

For the sphere $S$, presented as either of the two CCS's on the left of Figure \ref{fig:inv-exp1}, we have:
$$
\left \langle \emptyset  \left | Z_S \right |\emptyset  \right \rangle = \frac{1}{|G|}
$$
since in the extreme left surface, the set of colourings $\{ j\, |\, j\in G\}$ has cardinality $|G|$ and there are two internal vertices, and in the surface next to it, the set of colourings
$\{(j_1,j_2,k)\in G^3\, | \,k=1\}$ has cardinality $|G|^2$ and there are three internal vertices.

For the disk $D$ (third CCS from the left in Figure \ref{fig:inv-exp1}), we have:
$$
\left \langle \emptyset  \left | Z_D \right | g  \right \rangle = \frac{1}{|G|^{1/2}} D(g),
$$
where 
\be
D(g) :=\left\{ \begin{array}{cc} 1, & g=1 \\ 0, & g\neq 1 \end{array} \right . 
\label{eq:Ddef}
\ee
since the set of colourings $\{ k\in G\, |\, gkk^{-1}=g=1\}$ has cardinality $|G|$, if $g=1$, and $0$ otherwise, and there is one internal vertex and one external vertex.

For the cylinder $C$ (CCS on the right in Figure \ref{fig:inv-exp1}), noting that there are no internal vertices and two external vertices, we have
$$
\left \langle h  \left | Z_C \right | g  \right \rangle =  \frac{1}{|G|} C(g,h)
$$
where 
\be
C(g,h) := \# \{ k\in G\, |\, g= khk^{-1}\}\, .
\label{eq:Cdef}
\ee
For the torus $T$ (CCS on the left in Figure \ref{fig:inv-exp2}), we get:
$$
\left \langle \emptyset  \left | Z_T \right |\emptyset  \right \rangle = \frac{1}{|G|} \#\{(j_1,j_2)\in G^2\, | \, j_1j_2=j_2j_1\} \, .
$$
Finally, for the pants surface $P$ (CCS on the right in Figure \ref{fig:inv-exp2}) we get:
$$
\left \langle i,h  \left | Z_P \right | g  \right \rangle = \frac{1}{|G|^{3/2}} P(g,i,h),
$$
where
\be
P(g,i,h) := \# \{ (j_1,j_2)\in G^2\, |\, g= j_2ij_2^{-1}j_1hj_1^{-1}\}\, .
\label{eq:Pdef}
\ee
\begin{Remark}
We will be studying the three functions $D(g), \,C(g,h), \, P(g,i,h)$ more closely soon. Note that $C(g,h)$ tells us the ``extent'' to which $g$ and $h$ are conjugate, e.g. $C(g,h)=0$ if $g$ and $h$ are not conjugate, and $C(g,g)$ is the cardinality of the centraliser of $g$ in $G$. Likewise, comparing the expressions for $C(g,h)$ and $P(g,i,h)$, the latter function tells us the ``extent'' to which $g$ is 
``doubly conjugate'' to the pair of elements $i,h$.
\end{Remark}

\section{A gluing formula for the invariant and TQFT's for CCS's}
\label{sec:c}

\begin{figure}[htbp] 
\centerline{\relabelbox 
\epsfxsize 10cm
\epsfbox{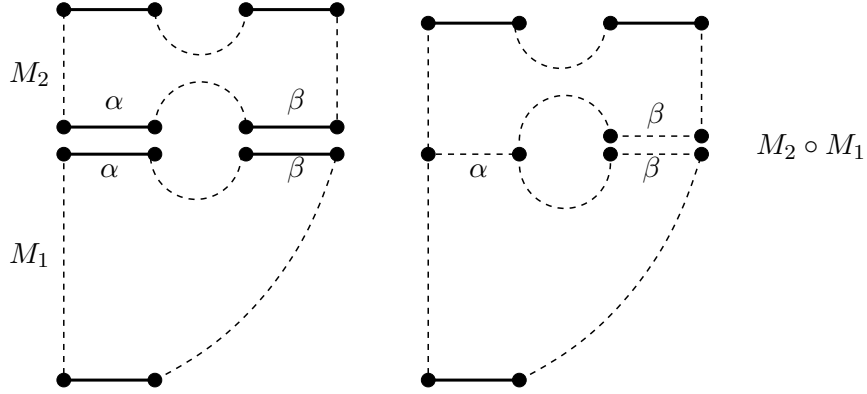}
\relabel{a}{$M_1$}
\relabel{b}{$M_2$}
\relabel{c}{$M_2\circ M_1$}
\relabel{h}{$\alpha$}
\relabel{i}{$\alpha$}
\relabel{e}{$\alpha$}
\relabel{j}{$\beta$}
\relabel{d}{$\beta$}
\relabel{f}{$\beta$}
\relabel{g}{$\beta$}
\endrelabelbox}
\caption{Gluing or composition of two CCS's}
\label{fig:gluing}
\end{figure}

When the outgoing boundary components of $M_1$ match the incoming boundary components of $M_2$, we can perform a gluing, or composition, to make a surface $M_2\circ M_1$. For the result to be a CCS, we adopt the following procedure, illustrated in Figure \ref{fig:gluing}. 
\vskip 0.1cm
\noindent {\em Gluing procedure:} we identify the shared boundary component furthest to the left (labelled $\alpha$ in the Figure), and the remaining shared boundary components (just one  in the Figure, labelled $\beta$ ) become cut 1-cells in the boundary of the planar representation of $M_2\circ M_1$.

\begin{Remark}
We could have obtained a different CCS by making $\beta$ internal instead of $\alpha$. In general, one could define a set of compositions between $M_1$ and $M_2$, specifying which shared boundary becomes internal. 
\end{Remark}

Suppose we have $M_1$ with $n$ incoming boundary components and $m>0$ outgoing boundary components, and $M_2$ with $m$ incoming boundary components and $p$ outgoing boundary components. Fixing the colourings of the ``in'' boundary components of $M_1$ and the ``out'' boundary components of $M_2$ the colourings of $M_2\circ M_1$ allow a priori any choice for the colourings of the $m$ intermediate 1-cell components. Thus we arrive at the following gluing formula for the invariant:
\begin{eqnarray}
\lefteqn{\left \langle i_1,\dots , i_p \left | Z_{M_2\circ M_1} \right |g_1,\dots , g_n   \right \rangle =} \nonumber \\  
%& & \sum_{\tiny \begin{array}{c}h_i\in G \\ i=1, \dots, m\end{array}}
& &  \sum_{\tiny h_i\in G }  
\left \langle i_1,\dots , i_p \left | Z_{M_2} \right | h_1,\dots , h_m  \right \rangle  
%\left \langle h_1,\dots , h_m \left | Z_{M_1} \right |g_1,\dots , g_n   \right \rangle \nonumber
%\\ 
\left \langle h_1,\dots , h_m \left | Z_{M_1} \right |g_1,\dots , g_n   \right \rangle 
\label{eq:gluing}
\end{eqnarray}
Note that the $m$ external vertices in the two factors on the right hand side each give rise to a factor $\frac{1}{|G|^{1/2}}$ in (\ref{def:Zgrp}), and these combine to give a factor $\frac{1}{|G|}$ for the corresponding internal vertex on the left hand side.

\vskip 0.5 cm

The gluing formula enables us to construct a natural TQFT, as already suggested by our notation for the invariants. We assign to each incoming or outgoing boundary of a CCS $M$, a vector space $V_i$ or $V_o$ over $\mathbb{R}$, whose basis consists of all $G$-colourings of the boundary components. The basis elements are written $|\,g_1,\dots , g_n  \left. \right \rangle$ or  
$\left \langle \right. h_1,\dots , h_m\,  |$, and the dimension of $V_i$ and $V_o$ is $|G|^n$ and $|G|^m$ respectively. To the CCS itself we assign the linear transformation $Z_M$ from $V_i$ to $V_o$, whose matrix elements with respect to these two bases are given by:
$$
 \left \langle h_1,\dots , h_m \left | Z_M \right |g_1,\dots , g_n   \right \rangle
$$
Thus from the gluing formula (\ref{eq:gluing}) we have the fundamental TQFT property:
\be
Z_{M_2\circ M_1} = Z_{M_2} \circ Z_{M_1}.
\label{tqft}
\ee
There is an important corollary of (\ref{tqft}), which expresses that the cylinder $C$ has assigned to it an idempotent, as should be expected from our general discussion of TQFT in Section 2 - see equation (\ref{eq:idemp}).
\begin{Proposition}
$$
Z_{C} \circ Z_{C} = Z_{C}.
$$
\label{ZC}
\end{Proposition}
\begin{Proof}
This holds due to the invariance of the matrix elements of $Z_C$, implying that they are the same for the simple cylinder as in Figure \ref{fig:inv-exp1}, and for the cylinder made of two cylinders glued together, which is obtained from the simple cylinder by a move of Type I followed by a move of Type II. Thus $Z_{C\circ C}=Z_C$ and the property follows. 
\end{Proof}

\begin{Remark} In terms of the function $C(g,h)$ defined in (\ref{eq:Cdef}), this corresponds to the equation:
\begin{equation}
\sum_{h \in G} C(g,h) C(h,i) = |G| \cdot C(g,i),
\label{eq:cyl}
\end{equation}
which we can also prove algebraically by straightforward means:

If $g$ and $i$ are not conjugate, then the right hand side of (\ref{eq:cyl}) is 0, and so is the left hand side, since for any $h$, either $C(g,h)=0$ or $C(h,i)=0$. It remains to consider the terms in (\ref{eq:cyl}) for which $g$, $h$ and $i$ are in the same conjugacy class. When $h$ is conjugate to $g$, i.e. there exists $m$ such that $h=mgm^{-1}$, we have $C(g,h)=C(g,g)$, since
$$
g=khk^{-1} \Leftrightarrow g= (km) g (km)^{-1}.
$$ 
Hence $C(g,h)=C(h,i)=C(g,i)$, and (\ref{eq:cyl}) reduces to 
\begin{equation}
\sum_{h \in G} C(g,h) = |G|,
\label{eq:sumC}
\end{equation}
which clearly holds, since, fixing $g$, every $k\in G$ belongs to precisely one set of the form $\left\{ k\in G : g=khk^{-1}  \right\}$.

If we denote the conjugacy class of $g$ in $G$ by $\overline{g}$, we then have an equation for the number of elements of $\overline{g}$:
\begin{equation}
\# \overline{g} = |G|/C(g,g),
\label{eq:card-g-bar}
\end{equation}
since $C(g,h) = C(g,g)$ for all $h\in \overline{g}$, and the number of non-zero terms in the sum on the l.h.s. of (\ref{eq:sumC})  is $\# \overline{g}$. 

Let ${\rm ConjClass}(G)$ denote the set of conjugacy classes of $G$. Then its cardinality is given by the following equation, which we will be using shortly:
\begin{equation}
\# {\rm ConjClass}(G) = \frac{1}{|G|^2} \sum_{g,h \in G} C(g,h)^2.
\label{eq:conj-cl-C}
\end{equation}
The double sum on the r.h.s. decomposes into double sums where $g,h$ both belong to the same conjugacy class $\overline{g}$. Restricting to these terms, the r.h.s. of (\ref{eq:conj-cl-C})
becomes:
$$
\frac{1}{|G|^2} \sum_{g,h \in \overline{g}} C(g,h)^2=
\frac{(\# \overline{g})^2 (C(g,g)^2}{|G|^2} = 1
$$
and collecting the contributions from each conjugacy class, we get equation (\ref{eq:conj-cl-C}).

\end{Remark}

\begin{Remark}
Note that this property of $C(g,h)$ distinguishes our TQFT's from the ordinary 2-dimensional TQFT's, which were shown by Abrams to be equivalent to Frobenius algebras \cite{at}. For these the linear mapping assigned to the cylinder is simply the identity map on the vector space associated to the circle. In our case it is clear that the invariant subspaces of $Z_C$ are those spanned by collections of elements of $G$ which constitute a conjugacy class. Within each invariant subspace there is just one eigenvector  of $Z_C$ with eigenvalue 1, namely $\sum_{g\in \cal C}\left | g \right \rangle $ where $\cal C$ denotes a conjugacy class of $G$.
\end{Remark}

\begin{figure}[htbp] 
\centerline{\relabelbox 
\epsfxsize 12cm
\epsfbox{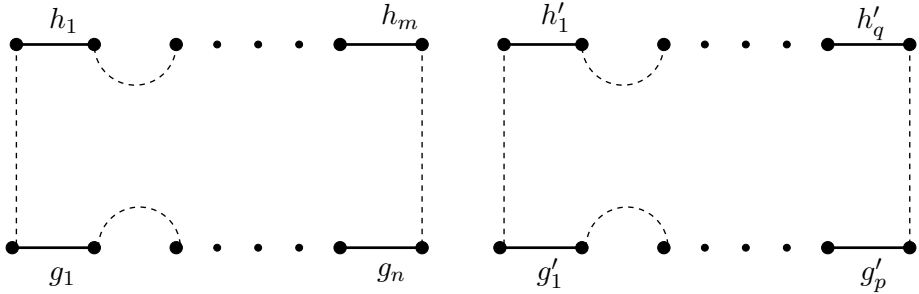}
\relabel{a}{$h_1$}
\relabel{b}{$h_m$}
\relabel{c}{$h_1'$}
\relabel{d}{$h_q'$}
\relabel{f}{$g_1$}
\relabel{g}{$g_n$}
\relabel{h}{$g_1'$}
\relabel{i}{$g_p'$}
\endrelabelbox}
\caption{The disjoint union of two CCS's}
\label{fig:disj-union}
\end{figure}

To complete the TQFT description, we also need to address the disjoint union of two CCS's $M_1$ and $M_2$, denoted $M_1  \sqcup M_2$, which is simply their juxtaposition, as in Figure \ref{fig:disj-union}. 

The number of $G$-colourings of   $M_1  \sqcup M_2$ is clearly the product of the number of $G$-colourings of $M_1$ and $M_2$ separately, fixing the colours of the boundary components. Defining the invariant for $M_1  \sqcup M_2$ in exactly the same way as for a single CCS $M$, we therefore have the property:

\begin{eqnarray}
\lefteqn{ \langle h_1,\cdots , h_{m}, h_1', \cdots , h_q' \left | Z_{M_1  \sqcup M_2} \right |g_1,\cdots , g_{n}, g_1', \cdots , g_p'    \rangle =} \nonumber \\  
& &    
 \langle h_1,\cdots , h_m \left | Z_{M_1} \right | g_1,\cdots , g_n   \rangle   \, . \,
 \langle h'_{1},\cdots , h'_{q} \left | Z_{M_2} \right |g'_{1},\cdots , g'_{p}    \rangle  
\label{eq:disj-union}
\end{eqnarray}

This multiplicative property can also be expressed more concisely using the tensor product of linear transformations:
$$
Z_{M_1  \sqcup M_2} = Z_{M_1 }\otimes  Z_{M_2} \, . 
$$
\vskip 0.2cm

Our main point, to be developed in the next section, is that the gluing formula (\ref{eq:gluing}) and TQFT relations (\ref{tqft}) and (\ref{eq:disj-union}) provide a powerful and general approach to proving algebraic properties of the functions $D(g)$, $C(g,h)$ and $P(g,i,h)$, like equation (\ref{eq:cyl}), by topological means.

\section{Properties of $D$, $C$ and $P$ from topology}
\label{sec:DCP}

\subsection{Initial properties of $D$, $C$ and $P$}
\label{ssec:DCPinitial}

The functions $D(g)$, $C(g,h)$ and $P(g,i,h)$ have an algebraic definition, but enjoy properties coming from their origin as invariants of CCS's, associated to the disk $D$, the cylinder $C$ and the pants surface $P$, respectively - see equations (\ref{eq:Ddef}),  (\ref{eq:Cdef}), (\ref{eq:Pdef}). In this final section we will analyse some relations that these functions obey, exploiting the invariance of $Z_M$ for suitable choices of $M$.

\begin{figure}[htbp] 
\centerline{\relabelbox 
\epsfysize 3.2cm
\epsfbox{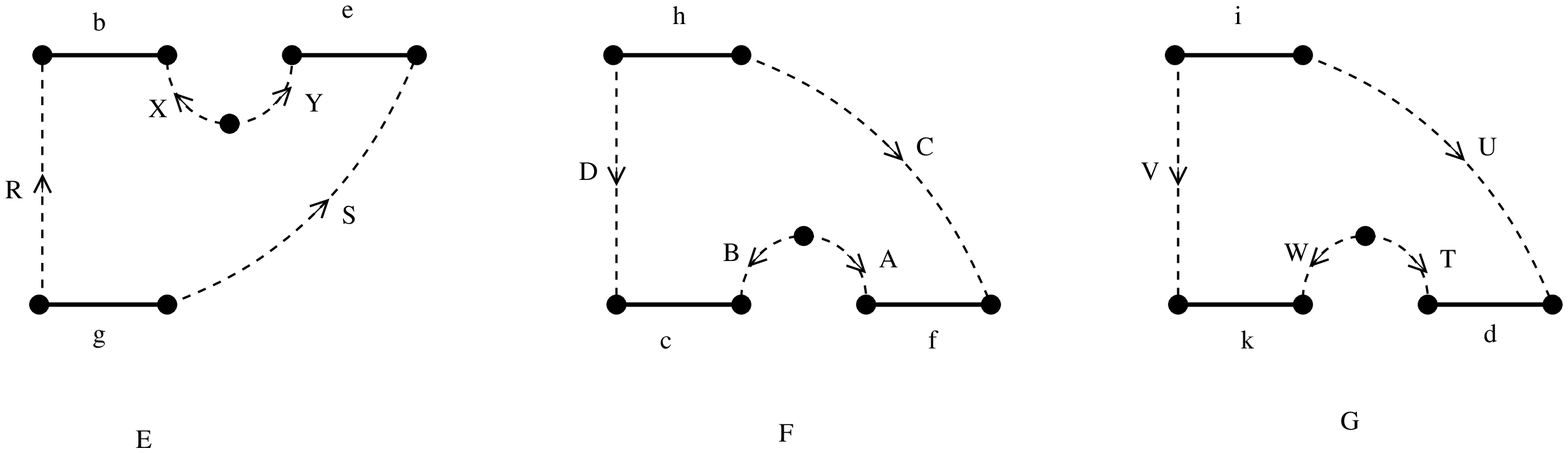}
\adjustrelabel <-2pt,0pt> {R}{$\alpha$}
\relabel{S}{$\beta$}
\adjustrelabel <-2pt,0pt> {X}{$\alpha$}
\relabel{Y}{$\beta$}
\adjustrelabel <-2pt,0pt> {D}{$\alpha$}
\adjustrelabel <-2pt,0pt> {B}{$\alpha$}
\relabel{A}{$\beta$}
\relabel{C}{$\beta$}
\relabel{T}{$\alpha$}
\relabel{U}{$\alpha$}
\adjustrelabel <-2pt,0pt> {V}{$\beta$}
\adjustrelabel <-2pt,0pt> {W}{$\beta$}
\relabel{g}{$\gamma$}
\relabel{h}{$\gamma$}
\relabel{i}{$\gamma$}
\relabel{b}{$\delta$}
\adjustrelabel <0pt,-1pt> {c}{$\delta$}
\relabel{d}{$\delta$}
\adjustrelabel <0pt,-1pt> {e}{$\epsilon$}
\relabel{f}{$\epsilon$}
\relabel{k}{$\epsilon$}
\relabel{E}{$P$}
\relabel{F}{$\overline{P}$}
\relabel{G}{$P_r$}
\endrelabelbox}
\caption{\label{fig:pants-ops} Reflection and 180 degree rotation of $P$}
\end{figure}

Because of the planar representation that a CCS possesses by definition, there are two natural symmetry operations on CCS's, namely reflection in a horizontal axis and rotation by 180 degrees around an axis perpendicular to the plane. Both of these operations exchange the incoming and outgoing boundary components, and the 180 degree rotation also inverts the orientation of the boundary components, so this has to be adjusted after the rotation. We denote the CCS obtained from $M$ by reflection and by 180 degree rotation, as $\overline{M}$ and $M_r$ respectively.  See Figure \ref{fig:pants-ops} for an example of these operations applied to the pants surface $P$.

The invariant $Z_M$ has the following properties under these operations.

\begin{Proposition}
$$
\left \langle g_1,\dots , g_n  \left | Z_{\overline{M}} \right |h_1,\dots , h_m    \right \rangle
=
\left \langle h_1,\dots , h_m \left | Z_M \right |g_1,\dots , g_n   \right \rangle \, ,
$$
$$
\left \langle g_n^{-1},\dots , g_1^{-1}  \left | Z_{M_r} \right |h_m^{-1},\dots , h_1^{-1}    \right \rangle
=
\left \langle h_1,\dots , h_m \left | Z_M \right |g_1,\dots , g_n   \right \rangle \, .
$$
\label{prop:pants-ops}
\end{Proposition}
\begin{Proof}
There is a one-to-one correspondence between the $G$-colourings of $M$ and those of $\overline{M}$ and $M_r$, as is clear by performing the symmetry operation on the whole labelled planar representation, and in the case of $M_r$, replacing $g_i$ by $g_i^{-1}$, $h_j$ by $h_j^{-1}$, and inverting the order of the respective boundary components.
\end{Proof}
\vskip 0.1cm

Since $D_r= \overline{D}$, $C_r= \overline{C} = C $ and $P_r= \overline{P}$, we deduce the following properties of the functions $D$, $C$ and $P$.

\begin{Corollary}
\begin{eqnarray}
D(g^{-1}) & = & D(g) \\
C(h^{-1},g^{-1})& = & C(h,g) \,\, = \,\, C(g,h)   \\
 P(g^{-1},h^{-1},i^{-1}) &= & P(g,i,h)
\end{eqnarray}
\end{Corollary}
\begin{Proof}
As an example, the last relation is derived as follows:
\begin{eqnarray*}
\frac{1}{|G|^{3/2}} P(g^{-1},h^{-1},i^{-1}) & = & \left \langle h^{-1},i^{-1}   \left | Z_{P} \right |  g^{-1} \right \rangle \\
& = & \left \langle g^{-1}  \left | Z_{\overline{P}} \right | h^{-1},i^{-1}   \right \rangle \\
&  = & \left \langle g^{-1}  \left | Z_{P_r} \right | h^{-1},i^{-1}   \right \rangle \\
&  = & \left \langle i,h  \left | Z_{{P}} \right | g   \right \rangle \\
&  = & \frac{1}{|G|^{3/2}}  P(g,i,h)
\end{eqnarray*}
\end{Proof}

\begin{figure}[htbp] 
\centerline{\relabelbox 
\epsfysize 2.2cm
\epsfbox{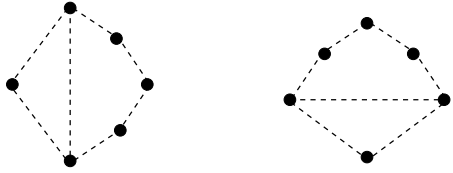}
\endrelabelbox}
\caption{Example of $M$ and $M'$ for Proposition \ref{prop:empty-rotn}}
\label{fig:empty-rotn}
\end{figure}

When $M$ has empty in and out boundaries, the invariance under rotations generalizes to the following property, which is proved analogously to the previous proposition (see Figure \ref{fig:empty-rotn})

\begin{Proposition}
Let $M$ and $M'$ be two CCS's with empty in and out boundaries, and sharing the same structure as a CCS apart from a different choice of incoming and/or outgoing 0-cell. Then:
$$
\left \langle \emptyset  \left | Z_{M'} \right | \emptyset   \right \rangle
=
\left \langle \emptyset  \left | Z_M \right | \emptyset  \right \rangle \, .
$$
\label{prop:empty-rotn}
\end{Proposition}

To complete this initial set of results, we have the following property, allowing a simplification of the presentation of a CCS, when a disk is glued into an open boundary component.

\begin{figure}[htbp] 
\centerline{\relabelbox 
\epsfysize 4cm
\epsfbox{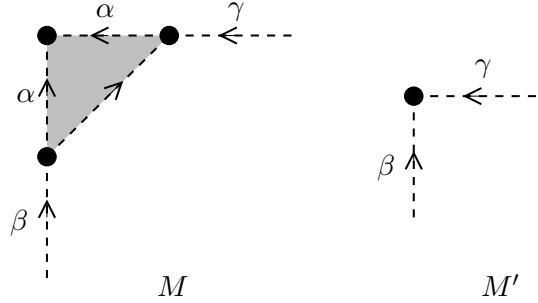}
\relabel{T}{$\alpha$}
\relabel{X}{$\beta$}
\relabel{U}{$\alpha$}
\relabel{a}{$\beta$}
\adjustrelabel <-1pt,1pt> {Y}{$\gamma$}
\relabel{b}{$\gamma$}
\relabel{Z}{$M$}
\relabel{W}{$M'$}
\endrelabelbox}
\caption{Illustration for Proposition \ref{prop:disk-out} }
\label{fig:disk-out}
\end{figure}

\begin{Proposition}
Suppose $M$ and $M'$ are two CCS's with the same structure as a CCS apart from one 2-cell in $M$ which is shrunk to a 0-cell in $M'$, like the shaded 2-cell in Figure \ref{fig:disk-out}. 
Then the invariants for $M$ and $M'$ are the same.
\label{prop:disk-out}
\end{Proposition}

\begin{Proof}
For any $G$-colouring of $M'$, there are $|G|$ compatible $G$-colourings of $M$, by choosing an arbitrary element of $G$ to assign to the 1-cell $\alpha$. However, this extra factor in the number of colourings is cancelled since $M$ has an extra internal 0-cell, which contributes a factor $1/|G|$ in the definition of $Z_M$. Two of the 0-cells in $M$ are the same, both being the starting point of 
$\alpha$, so that $M$ indeed  only has one extra internal 0-cell, not two, compared to $M'$.  Note: the orientations on $\alpha$, $\beta$, $\gamma$ chosen in Figure 
\ref{fig:disk-out} are just for illustration and may be inverted.
\end{Proof}

%%%%%%%%%%%%%%%%%%%%
%%%%%%%%%%%%%%%%%%%%
% Uma subseco

\subsection{An example}
\label{ssec:DCPtorusexp}

\begin{figure}[htbp] 
\centerline{\relabelbox 
\epsfysize 7cm
\epsfbox{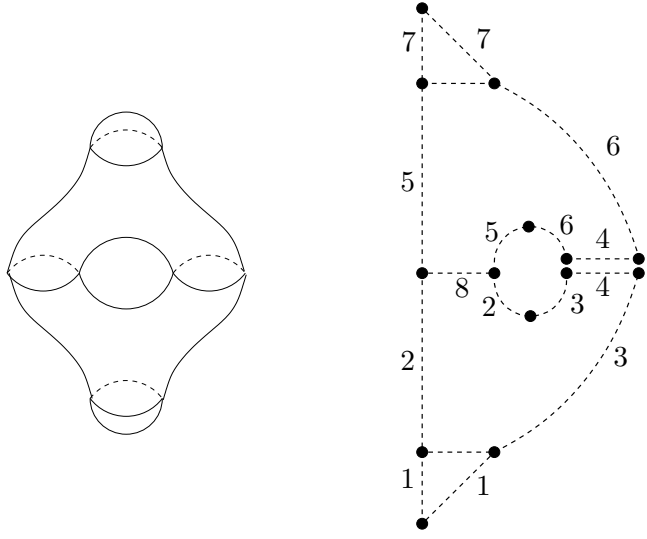}
\relabel{a}{$1$}
\relabel{b}{$1$}
\relabel{c}{$2$}
\adjustrelabel <-2pt,0pt> {d}{$2$}
\relabel{e}{$3$}
\relabel{h}{$3$}
\adjustrelabel <0pt,-2pt> {f}{$4$}
\relabel{g}{$4$}
\relabel{i}{$5$}
\adjustrelabel <-1pt,0pt> {j}{$5$}
\relabel{k}{$6$}
\relabel{l}{$6$}
\adjustrelabel <-1pt,0pt> {m}{$7$}
\relabel{n}{$7$}
\adjustrelabel <0pt,-1pt> {o}{$8$}
\endrelabelbox}
\caption{The torus obtained from gluing two pants surfaces and two disks }
\label{fig:torus-example}
\end{figure}

In this subsection we discuss the implications, in our context, of a simple topological fact, namely that the torus can be obtained by gluing two pairs of pants surfaces and two disks, as in Figure \ref{fig:torus-example} on the left. 

In terms of CCS's, the gluing of these four surfaces can be represented as in Figure \ref{fig:torus-example} on the right (we have omitted displaying the orientations so as not to clutter up the figure). The corresponding topological relation can be expressed in the language of TQFT, as:
$$
\left \langle \emptyset  \left | Z_{D} \circ Z_{\overline{P}} \circ  Z_{P} \circ Z_{\overline{D}} \right | \emptyset   \right \rangle
=
\left \langle \emptyset  \left | Z_T \right | \emptyset  \right \rangle \, .
$$

This relation may be proved both algebraically, by direct calculations, or by using properties of the invariant. First we give the algebraic derivation.

\begin{eqnarray}
l.h.s. & = & \frac{1}{|G|}  \sum_{h,k \in G} \frac{1}{|G|}  \left \langle 1 \left | Z_P \right | h, k  \right \rangle   \left \langle h, k  \left | Z_P \right | 1  \right \rangle   \nonumber  \\
& = &\frac{1}{|G|^4} \sum_{h,k \in G}  \# \left \{ j_1, j_2 \in G : j_1 h^{-1} j_1^{-1} j_2 k^{-1} j_2^{-1} = 1\right \} ^2  \nonumber  \\
& = & \frac{1}{|G|^4} \sum_{h,k \in G} \# \left \{ j_1, j_2 \in G : h^{-1} = j_1^{-1} j_2 k j_2^{-1} j_1 \right \} ^2  \nonumber  \\
& = & \frac{1}{|G|^4} \sum_{h,k \in G} \left [ |G|\# \left \{ j \in G : h^{-1} = j k j^{-1} \right \} \right ] ^2 \nonumber  \\ 
& = & \frac{1}{|G|^2} \sum_{h,k \in G} C(h^{-1},k) \cdot C(h^{-1},k) \nonumber  \\ 
& = & \frac{1}{|G|^2} \sum_{k \in G} |G| C(k,k) \nonumber \\ 
& = & \frac{1}{|G|} \sum_{k \in G} \# \left \{ j \in G : k j = j k \right \} \nonumber \\ 
& = & \frac{1}{|G|} \# \left \{ j, k \in G : k j = j k \right \} = r.h.s.
\label{T=2P+2D}
\end{eqnarray}

\begin{figure}[htbp] 
\centerline{\relabelbox 
\epsfysize 6cm
\epsfbox{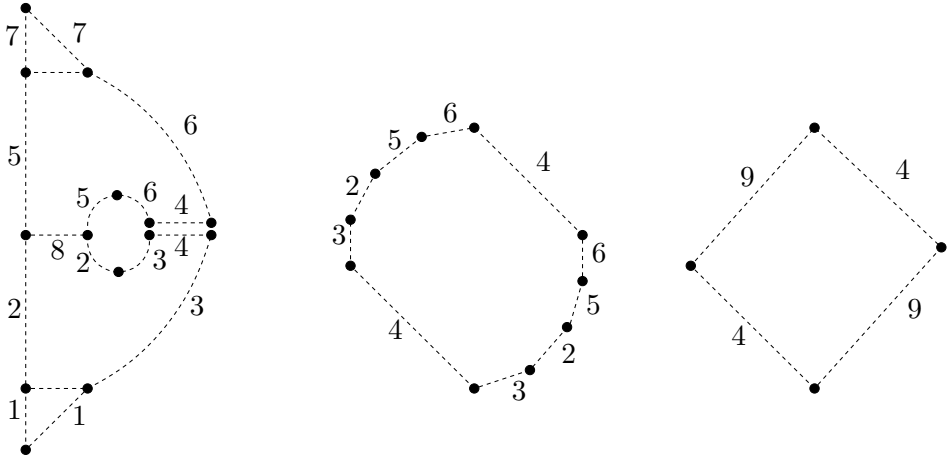}
\relabel{a}{$1$}
\relabel{b}{$1$}
\relabel{c}{$2$}
\adjustrelabel <-2pt,0pt> {d}{$2$}
\relabel{e}{$3$}
\relabel{h}{$3$}
\adjustrelabel <0pt,-2pt> {f}{$4$}
\relabel{g}{$4$}
\relabel{i}{$5$}
\adjustrelabel <-2pt,0pt> {j}{$5$}
\relabel{k}{$6$}
\relabel{l}{$6$}
\adjustrelabel <-2pt,0pt> {m}{$7$}
\relabel{n}{$7$}
\adjustrelabel <-1pt,-1pt> {o}{$4$}
\relabel{p}{$3$}
\adjustrelabel <-2pt,0pt> {q}{$2$}
\relabel{r}{$5$}
\relabel{s}{$6$}
\relabel{t}{$4$}
\adjustrelabel <0pt,-1pt> {u}{$6$}
\adjustrelabel <-1pt,-1pt> {v}{$5$}
\adjustrelabel <-1pt,-3pt> {w}{$2$}
\adjustrelabel <-1pt,-2pt> {x}{$3$}
\adjustrelabel <-1pt,-1pt> {y}{$4$}
\relabel{z}{$9$}
\adjustrelabel <-1pt,-1pt> {A}{$4$}
\relabel{B}{$9$}
\adjustrelabel <0pt,-1pt> {C}{$8$}
\endrelabelbox}
\caption{Moves on CCS representations}
\label{fig:CCS-torus-moves}
\end{figure}

Alternatively we display in Figure \ref{fig:CCS-torus-moves}, a sequence of moves on CCS representations showing the equality. 

In the first step we use Move II and  Prop. \ref{prop:disk-out} to remove 1-cells 1,7 and 8, and Prop. \ref{prop:empty-rotn} to rearrange the ``in'' and ``out'' 0-cells. Then we use Move I to replace four 1-cells with a single 1-cell.

Irrespective of the method of proof, we obtain an identity for the functions $D$ and $P$ in terms of $|G|$ and the invariant for the torus: 

$$
 \sum_{g,h,k,l \in G} D(g)P(g,h,k)P(l,h,k) D(l) = |G|^3 \# \left \{ j, k \in G : k j = j k \right \} 
$$

Although the function $C$ does not enter explicitly in this identity, the algebraic derivation shows that it appears along the way. We will next show that $C$ in fact depends on $D$ and $P$, i.e. is not an independent function.

\begin{figure}[htbp] 
\centerline{\relabelbox 
\epsfysize 4.2cm
\epsfbox{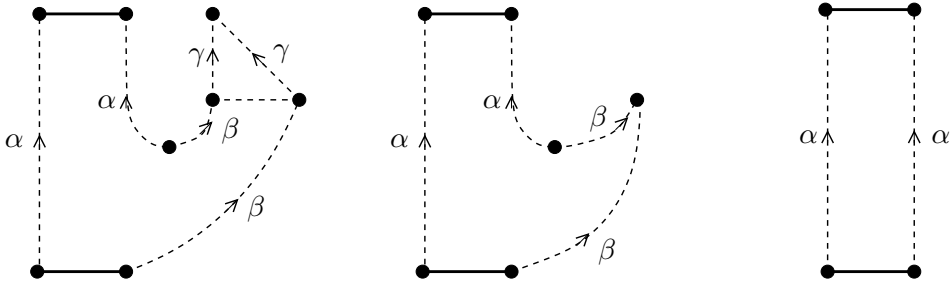}
\relabel{R}{$\alpha$}
\relabel{S}{$\beta$}
\adjustrelabel <-1pt,0pt> {X}{$\alpha$}
\relabel{Y}{$\beta$}
\adjustrelabel <-1pt,0pt> {e}{$\gamma$}
\relabel{h}{$\gamma$}
\relabel{D}{$\alpha$}
\adjustrelabel <-1pt,0pt> {B}{$\alpha$}
\adjustrelabel <-1pt,0pt> {V}{$\alpha$}
\relabel{U}{$\alpha$}
\relabel{A}{$\beta$}
\relabel{C}{$\beta$}
\endrelabelbox}
\caption{The cylinder from the pants surface and the disk}
\label{fig:C=PD}
\end{figure}

It is enough to consider the CCS $M$ displayed on the left of Figure \ref{fig:C=PD}, and the moves in the same Figure using Proposition \ref{prop:disk-out} and Move II, showing that its invariant $Z_M$ is the same as that of the cylinder $Z_C$. 

Thus for all $g,i\in G$
$$
\left \langle i  \left | Z_M \right | g  \right \rangle =
\left \langle i \left | Z_C \right | g  \right \rangle =
\frac{1}{|G|} C(g,i)
$$
The number of $|G|$-colourings of $M$ can be expressed as:
$$
\# {\cal C}_M (g,i) =  \sum_{h \in G} P(g,i,h) |G| D(h) =
|G| P(g,i,1)
$$ 
Since 
$$
\left \langle i  \left | Z_M \right | g  \right \rangle =
\frac{1}{|G|^3}\# {\cal C}_M (g,i)
$$
we arrive at the formula:
$$
C(g,i)= \frac{1}{|G|} \sum_{h \in G} P(g,i,h)  D(h) = 
\frac{1}{|G|} P(g,i,1).
$$

We conclude this subsection with an interesting result related to the torus. Returning to halfway through the derivation of (\ref{T=2P+2D}) and relabelling, we obtain:
\begin{equation}
\frac{1}{|G|} \# \left \{ j, k \in G : k j = j k \right \} =
\frac{1}{|G|^2} \sum_{g,h \in G} (C(g,h))^2
\label{eq:T=2C}
\end{equation}
which reflects the topological fact that the torus is obtained by gluing two cylinders together.

Now, for a finite group $G$, its commuting fraction is defined to be the quotient:
$$
\frac{\left \{ j, k \in G : k j = j k \right \}}{|G|^2}
$$
i.e. the number of commuting pairs of elements of $G$ divided by the overall number of pairs. Combining equation (\ref{eq:T=2C}) and 
equation (\ref{eq:conj-cl-C}), we have a topological proof for the following proposition, which is well-known to group theorists.

\begin{Proposition}
The number of conjugacy classes of $G$ is equal to the commuting fraction of $G$ times the order of $G$. 
\label{cf-prop}
\end{Proposition}
We are grateful to an anonymous referee for drawing our attention to the commuting fraction and to this property.

\subsection{The invariant and pants decompositions}
\label{ssec:inv+pants}

\begin{figure}[htbp] 
\centerline{\relabelbox 
\epsfxsize 12cm
\epsfbox{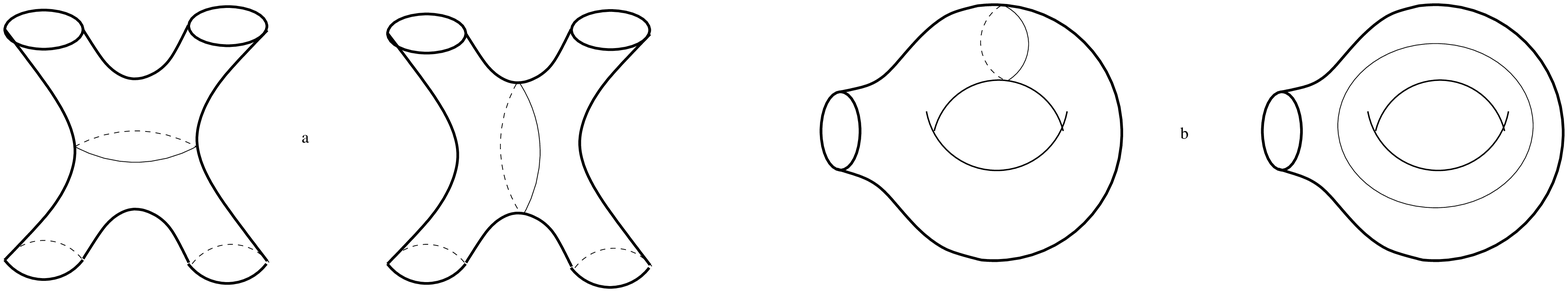}
\relabel{a}{$\leftrightarrow$}
\relabel{b}{$\leftrightarrow$}
\endrelabelbox}
\caption{Local moves on pants decompositions}
\label{fig:pants-moves}
\end{figure}

Apart from a small number of simple cases, all surfaces with boundary can be realized by gluing together a number of pants surfaces (the exceptions are the sphere, disk, cylinder and torus). A result due to Hatcher and Thurston \cite{ht} states that any two pants decompositions of the same surface can be related by a finite sequence of two types of local moves, depicted in Figure \ref{fig:pants-moves}.

The surface in the second move is a torus with one puncture, obtained by identifying two boundary components of a pants surface in two different ways.

In this subsection we show how these moves are described in the context of CCS's, and how the invariance of $Z_M$ gives rise to equations for the function $P$, associated with the pants surface.

\begin{figure}[htbp] 
\centerline{\relabelbox 
\epsfysize 5cm
\epsfbox{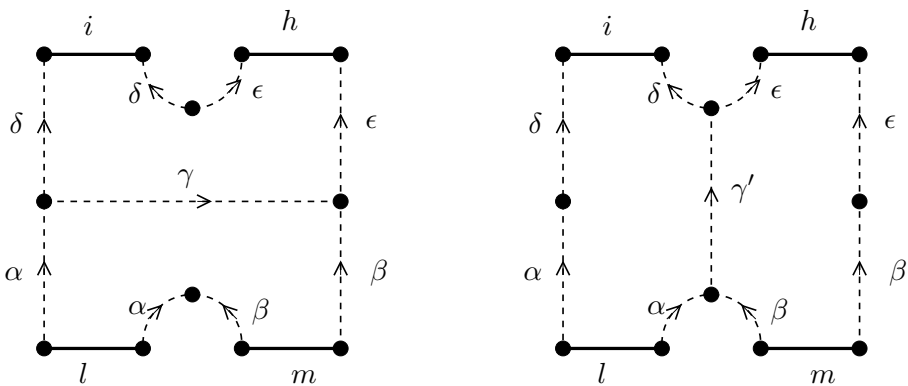}
\relabel{b}{$i$}
\relabel{c}{$l$}
\adjustrelabel <-0pt,-2pt> {e}{$h$}
\relabel{f}{$m$}
\relabel{D}{$\alpha$}
\relabel{B}{$\alpha$}
\relabel{A}{$\beta$}
\relabel{C}{$\beta$}
\relabel{R}{$\delta$}
\relabel{X}{$\delta$}
\relabel{Y}{$\epsilon$}
\relabel{S}{$\epsilon$}
\relabel{m}{$i$}
\relabel{n}{$l$}
\adjustrelabel <0pt,-2pt> {o}{$h$}
\relabel{p}{$m$}
\relabel{M}{$\alpha$}
\relabel{N}{$\alpha$}
\relabel{P}{$\beta$}
\relabel{Q}{$\beta$}
\relabel{I}{$\delta$}
\relabel{J}{$\delta$}
\relabel{K}{$\epsilon$}
\relabel{L}{$\epsilon$}
\relabel{h}{$\gamma$}
\relabel{i}{$\gamma'$}
\endrelabelbox}
\caption{The first pants move using CCS's}
\label{fig:2P-move}
\end{figure}

Starting with the first move, we represent it as a move on CCS's in Figure \ref{fig:2P-move}.

\begin{figure}[htbp] 
\centerline{\relabelbox 
\epsfysize 3cm
\epsfbox{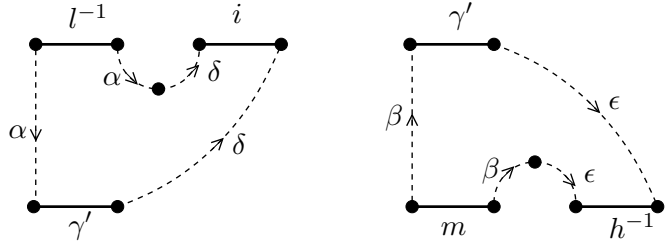}
\adjustrelabel <-1pt,0pt> {R}{$\alpha$}
\relabel{S}{$\delta$}
\adjustrelabel <-1pt,0pt> {X}{$\alpha$}
\relabel{Y}{$\delta$}
\relabel{A}{$\epsilon$}
\relabel{C}{$\epsilon$}
\adjustrelabel <-1pt,0pt> {D}{$\beta$}
\relabel{B}{$\beta$}
\adjustrelabel <0pt,-1pt> {g}{$\gamma'$}
\relabel{h}{$\gamma'$}
\relabel{b}{$l^{-1}$}
\adjustrelabel <0pt,-1pt> {e}{$i$}
\relabel{c}{$m$}
\adjustrelabel <-3pt,-1pt> {f}{$h^{-1}$}
\endrelabelbox}
\caption{The pants surfaces for the right-hand CCS in Figure \ref{fig:2P-move}}
\label{fig:2P-move-ill}
\end{figure}

These two CCS's are clearly related by two Type II moves; removing the 1-cell labelled $\gamma$ and inserting the 1-cell labelled $\gamma'$, or vice-versa. Thus the invariant is the same for both. The CCS on the left is a straightforward composition of surfaces: $P\circ \overline{P}$. The CCS on the right we can also recognize to be two pants surfaces glued together along a 1-cell, as illustrated in 
Figure \ref{fig:2P-move-ill}.

\noindent Here we have adjusted the group labels of two boundary components, because of the change of orientation. Thus we get the following identity for $P$, which holds for all $i,h,l,m \in G$:
\begin{equation}
\sum_{g\in G}  P(g,i,h) P(g,l,m) = \sum_{g\in G}   P(g,l^{-1},i) P(g,m,h^{-1})
\label{eq:2P-formula} 
\end{equation}

\begin{figure}[htbp] 
\centerline{\relabelbox 
\epsfysize 3cm
\epsfbox{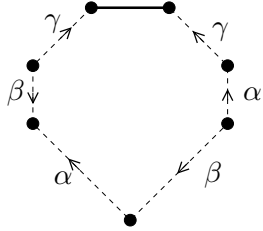}
\adjustrelabel <-1pt,-1pt> {R}{$\alpha$}
\relabel{X}{$\alpha$}
\adjustrelabel <-1pt,0pt> {D}{$\beta$}
\relabel{B}{$\beta$}
\relabel{g}{$\gamma$}
\relabel{h}{$\gamma$}
\endrelabelbox}
\caption{The torus with one puncture as a CCS}
\label{fig:1P-move-ill}
\end{figure}

For the second move we represent the torus with one disk removed as a CCS, $M$, in Figure \ref{fig:1P-move-ill}.

\begin{figure}[htbp] 
\centerline{\relabelbox 
\epsfysize 3cm
\epsfbox{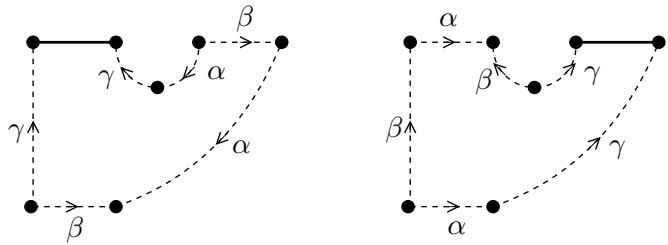}
\relabel{R}{$\gamma$}
\relabel{S}{$\alpha$}
\adjustrelabel <-2pt,0pt> {X}{$\gamma$}
\relabel{Y}{$\alpha$}
\adjustrelabel <0pt,-2pt> {g}{$\beta$}
\relabel{e}{$\beta$}
\relabel{A}{$\gamma$}
\relabel{C}{$\gamma$}
\relabel{D}{$\beta$}
\adjustrelabel <-2pt,-2pt> {B}{$\beta$}
\adjustrelabel <0pt,2pt> {h}{$\alpha$}
\relabel{c}{$\alpha$}
\endrelabelbox}
\caption{Illustration for the second pants move}
\label{fig:1P-move}
\end{figure}

Indeed if we glue a disk to the boundary component of $M$, we obtain a CCS which clearly reduces to the torus after applying Proposition \ref{prop:disk-out} twice. We now represent $M$ in two different ways as coming from a pants surface with two boundary components identified - see Figure \ref{fig:1P-move}.

Thus the invariant $\left \langle h  \left | Z_M \right | \emptyset  \right \rangle =
\frac{1}{|G|^{3/2}}\# {\cal C}_M (h)$ can be expressed in two different ways in terms of the function $P$, giving rise to the identity:

\begin{equation}
\sum_{g\in G}  P(g,h,g)  = \sum_{g\in G}   P(g,g,h)
\label{eq:1P-formula} 
\end{equation}

\section{Conclusions and Final Remarks}

Viewing our results from the group theory perspective, we have been led to introduce two integer-valued functions $C$ and $P$, which depend on two or three $G$-elements, respectively. The function $C$ tells us, first of all, whether its two arguments are conjugate or not, i.e. $C(g,h)\neq 0$ or $C(g,h)= 0$, but also gives a measure of conjugacy between the elements $g$ and $h$ by counting the number of elements $k$ such that $g=khk^{-1}$. Likewise the values $P(g,i,h)$ of the function $P$ tell us whether or not the element $g$ is in a ``doubly conjugate'' relation with the pair of elements $i$ and $h$, i.e. if there exist $j_1,\,j_2$ such that $g=j_1ij_1^{-1}j_2hj_2^{-1}$, and if so (i.e. $P(g,i,h)\neq 0$), this value tells us the extent to which this double conjugacy holds, by counting the possibilities for $j_1$ and $j_2$. We have derived properties of the functions $C$ and $P$ by using topological reasoning.

From the topology perspective, we have developed a general way of regarding surfaces as cut up into a planar shape, which may be a nice viewpoint in a wider sense. For instance, in the context of Stokes' theorem for flux integrals, the image of a local parametrization of a 2-dimensional manifold  with boundary $M$ is the interior of a 2-cell of a CCS representation of $M$. 

Using $G$-colourings of the 1-cells, where $G$ is a finite group, 
 we have found invariants for our class of surfaces with boundary, and a TQFT setting for them which is not the same as the usual description of 2-dimensional TQFT for surfaces. It would be interesting to look for other constructions of invariants, or generalizations such as including non-orientable surfaces. One construction which we are preparing \cite{bp2} is to use colourings by elements of a finite 2-group (i.e. a crossed module of finite groups) in order to obtain the invariants. In this case, there are two different finite groups $G$ and $H$, and the 1-cells of a CCS $M$ get elements of $G$ assigned to them, whilst the 2-cells of $M$ get elements of $H$ assigned to them, subject to a generalization of the flatness condition.

\section{Acknowledgments}
We are grateful to an anonymous referee for drawing our attention to the commuting fraction of a group $G$ and the possibility of using a topological approach to prove Proposition \ref{cf-prop}.

This article and its sequel \cite{bp2} are based on a research project carried out by Diogo Bragança under the supervision of  Roger Picken. The authors are grateful to the {\em Fundação Calouste Gulbenkian} for supporting this project through the programme {\em Novos Talentos em Matemática}, which aims to stimulate undergraduate research in mathematics.  

Diogo Bragança would like to thank his family for all the support they gave him, and his supervisor for guiding him throughout the project. 

Roger Picken is grateful to Dan Christensen and Jeffrey Morton for useful discussions and suggestions. This work was funded in part by the Center for Mathematical Analysis, Geometry and Dynamical Systems (CAMGSD),
Instituto Superior T\'ecnico, Universidade de Lisboa, through the project UID/MAT/04459/2013, and by the {\em Fundação para a Ciência e a Tecnologia} (FCT, Portugal) through the
``Geometry and Mathematical Physics Project'', FCT EXCL/MAT-GEO/0222/2012.

\end{document}